\documentclass{article}
\usepackage{graphicx}
\usepackage{amssymb,amsmath}
\def\ac#1{{a({#1})}}

\def\timehm{\count31=\time \count32=\count31 \divide\count31 by 60
\number\count31 \multiply\count31 by 60 \advance\count32 by
-\count31 :\ifnum\count32<10 0\fi \number\count32}

\def\ep{\varepsilon}

\def\ap{\alpha}
\def\A{\Cal A}

\def\tm#1{\left(\text{mod }#1\right)}

\newtheorem{thm}{Theorem}
\newtheorem{cor}[thm]{Corollary}
\newtheorem{lem}[thm]{Lemma}

\newtheorem{prop}[thm]{Proposition}

\newcommand{\Gal}{{\rm Gal}}
\newcommand{\Frob}{{\rm Frob }}

\newcommand{\GL}{{\rm GL}}

\def\mod{\mathop{\rm mod}\nolimits}

\def\Q{\mathbb{Q}}

\def\A{{\cal A}}

\begin{document}

\title{{\bf Modular forms with large coefficient fields via congruences
}}
\author{Luis V. Dieulefait, Jorge Jim\'{e}nez-Urroz, Kenneth A. Ribet\footnote{first named author partially supported by  grant
MTM2009-07024 and by an ICREA Academia prize, second by DGICYT Grant
MTM2009-11068}
\\
 }
\date{\empty}

\maketitle

\vskip -20mm
%

%

\section{Introduction}
In this paper we will exploit the theory of congruences between
modular forms to deduce the existence of newforms (in particular,
cuspidal Hecke eigenforms)
 with levels of certain specific types having arbitrarily large coefficient fields. We will
 only consider
 newforms
 of weight $2$ and trivial nebentypus.\\
 If the level is allowed
  to be divisible by a large power $n$ of a fixed prime, or by the cube of a large prime $p$, then the
   coefficient fields of {\bf all} newforms of this level will grow with $n$ (with $p$, respectively) due to results of
    Hiroshi Saito (cf. \cite{Sa}, Corollary 3.4;  see also \cite{Br}) showing that
     the maximal real subfield of certain cyclotomic field whose degree grows
    with $n$ (with $p$, respectively) will be
     contained in these fields of coefficients. Thus, it is natural to deal with  the question when the levels
     are square-free or almost-square-free, i.e.,  square-free except for the
     fact that they are divisible by a fixed power of a small prime.\\

In the square-free case, for any given number $t$, we will prove that in levels which are the
product of exactly $t$ primes there are newforms with arbitrarily large coefficient fields. We
 will recall results of Mazur on reducible primes for newforms of prime level that give the
 case of $t=1$. Then, a generalization of these results will allow us to deduce  the case $t=2$.
For $t \geq 3$ we follow a completely different approach, namely, we
exploit congruences involving certain elliptic curves whose
 construction is on the one hand related  to Chen's celebrated results on (a partial answer to) Goldbach's conjecture (cf.
 \cite{Ch}) and  on the other hand inspired by Frey curves as in the proof of Fermat's
 Last Theorem (the diophantine problem that we will consider will be
 a sort of Fermat-Goldbach mixed problem).
 It is via the level lowering results in \cite{Ri} that the desired congruence will be guaranteed. The precise statement of our first
  main result is the following:\\

\begin{thm}
\label{squarefree} Let $B$ and $t$ be two given positive integers.
Then, there exist $t$ different primes $p_1, p_2, ...., p_t$ such
that if we call $N$ their product, in the space of cuspforms of
weight $2$, level $N$ and trivial nebentypus there exists a newform
$f$ whose field of coefficients $\Q_f$ satisfies:
$$ [\Q_f : \Q]  > B. $$
\end{thm}

In the almost-square-free case, we will consider levels $N=2^kp_1\dots p_t$ which are square-free except that they are divisible by a small power of $2$.
We will prove that for any fixed $t$, among newforms with
such levels  the fields of coefficients have unbounded degree. We
will use congruences with certain $Q$-curves constructed from
solutions to the problem of finding prime values attained by the
expression $(x^4+y^2)/c$. Again, these $Q$-curves will also have
some features inspired by Frey curves, and the existence of the
desired congruences will be a consequence of level lowering. For
$c=1$ it is a celebrated result of Friedlander and Iwaniec that
infinitely many primes are of the form  $x^4 + y^2$ (cf.
\cite{fi2}). Here we give the following generalization

\begin{thm}\label{mainintro} Let $B>0$ fixed, and $\Lambda$ the usual Von Mangold function.
 Then, we have uniformly in $c\le (\log x)^B$
\begin{equation}
\sum\sum_{\kern-10pt c|a^2+b^4\le cx}\Lambda((a^2+b^4)/c)=K(c)x^{3/4}+o\left(x^{3/4}\right)
\end{equation}
where $a,b$ run over positive integers,  and $K$ is completely explicit in terms of $c$.
\end{thm}
In Theorem \ref{main} in section 6 we  give the precise value of $K(c)$.
For our application to congruences between modular
forms we only need a mild version of the particular case with fixed $c$ of the form $c =5^\ell$.
 The following is a direct consequence of the previous theorem
 \begin{cor}
 \label{primevalues}
 Let $c$ be a positive odd integer.  Then there are infinitely many primes of the form $(x^4+y^2)/c$ if and only if
 $c$ can be written as the sum of two squares.\\
\end{cor}
The precise statement of our second main result, the one covering
the almost-square-free level case, is the following:
\begin{thm}
\label{almostsquarefree} Let $B$ and $t$ be two given positive
integers. Then, there exist $\alpha \in \{  5, 8\}$ and $t$
different odd primes $p_1, p_2, ...., p_t$  such that if we call $N$
the product of these $t$ primes, in the space of cuspforms of weight
$2$, level $2^\alpha N$ and trivial nebentypus there exists a
newform $f$ with field of coefficients $\Q_f$ satisfying:
$$ [\Q_f : \Q]  > B. $$
\end{thm}
Let us stress that  the results on prime values of  $(x^4+y^2)/c$,
 Theorem \ref{mainintro} and its corollary, are interesting in its own right,
 independently of the application  to finding newforms with large coefficient fields.

  \section{Theorem \ref{squarefree} for the case of prime level: Mazur's argument}

  Suppose that the level $N$ is prime and that $\ell >3$ is a prime that
divides $N - 1$. Then it is proved in \cite{Ma} that the prime is
Eisenstein, meaning that there is a newform $f$ of weight $2$ and
level $N$  such that if we call $K_f$ its field of coefficients
there is a prime $\lambda$ dividing $\ell$ in the ring of integers
of $K_f$ for which we have $a_p \equiv 1 + p \; \mod{\lambda}$ for
all primes $p$. The residual mod $\lambda$ Galois representation
attached to $f$ is reducible. In particular, we have
$$a_2 \equiv 3
\; \mod{\lambda} \qquad \quad (1)$$ The coefficients of the modular
form $f$ and those of any Galois conjugate $f^\sigma$ all satisfy
the bound $| a_p | \leq 2 \sqrt{p}$, in particular $a_2$ and all its
Galois conjugates have absolute value bounded above by $2 \sqrt{2}
 < 3$. Then, $a_2 -3$ is a non-zero algebraic integer, whose norm is
divisible by $\ell$ because of  congruence (1) and with absolute
value at most $(3 + 2 \sqrt{2})^{\deg k_f}$. Hence
$$ \ell \leq (3 + 2 \sqrt{2})^{\deg k_f}$$
Thus, $\deg K_f$ is bigger than a fixed constant times $\log \ell$.
Taking $\ell$ big and using Dirichlet's theorem to find an $ N
\equiv 1 \; \mod{\ell}$, we can make $\deg K_f$ as big as we
like. This proves the case of prime level ($t=1$) of Theorem \ref{squarefree}.\\

\section{A Frey curve adapted to Chen results, and the case  $t
\geq 3$ of Theorem \ref{squarefree}}

Let $\ell$ be a (large) prime number, and assuming for the moment
the truth of Goldbach's conjecture let us  write the even number
$2^{\ell + 4}$ as the sum of two prime numbers: $2^{\ell + 4} = p +
q$. Since $p$ and $q$ are clearly non-congruent modulo $4$, we
assume without loss of generality that $p \equiv 3 \; \mod{4}$. Let
$F$ be the semistable Frey curve associated to the triple $p , q,
2^{\ell + 4}$:
           $$ y^2 = x (x- p) (x+ 2^{\ell + 4})  $$

Its conductor is $2pq$, while its minimal discriminant is $\Delta =
(2^{\ell
+ 4} pq)^2/ 2^8 = (2^\ell pq)^2$.\\
The modularity of all semistable elliptic curves, proved by Wiles in \cite{Wi}, implies that
there is a newform $f$ of weight $2$ and level $2pq$ corresponding to $F$.\\
 The mod $\ell$ Galois representation $F[\ell]$ of
$\Gal(\bar{\Q}/\Q)$ is irreducible by results of Mazur, and
unramified at $2$ because the $2$-adic valuation of the discriminant
is divisible by $\ell$ (as in the original Frey curves related to
solutions of Fermat's Last Theorem). Although it comes initially from
a newform  $f$ of level $2pq$, by level-lowering (cf. \cite{Ri}) it arises also
from a newform $f'$ of level $pq$. \\
The trace of the action of $\Frob \; 2$ on $F[\ell]$ is $\pm (1 +
2)$, because this is the well-known necessary condition for
level-raising, i.e., for the existence of an $\ell$-adic Galois
representation with semistable ramification at $2$ providing a lift
of $F[\ell]$, and we have such a lift by construction: it is given by the Galois
action on the full $\ell$-adic Tate module $T_\ell(F)$ of the curve
$F$. So if we call $\{ a_p\}$ the coefficients of $f'$ we get
$$a_2 \equiv \pm 3 \; \mod{\lambda} \qquad \; (2)$$
for a prime $\lambda$
dividing $\ell$ in the field of coefficients of $f'$.
 From this congruence we
can argue as we did in the previous section using congruence (1) an
conclude easily that just by choosing the prime $\ell$ sufficiently
large we can construct newforms of weight $2$ and level $N= pq$ a
product of two
different primes with arbitrarily large field of coefficients.\\
Since Goldbach's conjecture remains open, in order to get an
unconditional result we need to move to the case of three primes in
the level. Using the results of Chen on Goldbach's problem (cf.
\cite{Ch}) we know that for $\ell$ sufficiently large $2^{\ell + 4}$
can be written as the sum of a prime and a ``pseudo-prime", i.e., a
number that is either a prime or the product of two different
primes. Then, in particular, one of the following is true for
infinitely many $\ell$: $2^{\ell + 4}$ can be written as the sum of
two primes $p$ and $q$, or $2^{\ell + 4}$ can be written as the sum
of a prime $p$ and the product of two primes $qr$. If the first is
true, as we have just seen, this will prove Theorem \ref{squarefree}
for levels which are the product of two primes, and if the second is
true a similar argument with the triple $p, qr, 2^{\ell + 4}$ shows
that the theorem is true for levels which are the product of three
primes.\\
Thus, to finish the proof of the case $t=3$ of Theorem
\ref{squarefree}, it remains to show that if the result holds for
$t=2$ then it also holds for $t=3$. But this is just an application
of raising the level (cf. \cite{Rilr}), because whenever we have a
modular form of level $pq$ and an irreducible mod $\ell$ Galois
representation attached to it, the same residual representation is
also realized in some newform of level $pqr$, as long as the prime
$r$ satisfies the required condition for level raising (and it is
well known that there are infinitely many primes $r$ that satisfy
this condition, cf. \cite{Rilr}). Thus, whenever we have found a
newform $f'$ as in the previous argument, of level $pq$ and
satisfying (2), there are also newforms with levels of the form
$pqr$ also satisfying (2) and from this the proof of the theorem for
the case $t=3$ follows exactly as explained
above.\\
To treat the case of more than three primes, we modify the argument
above by further raising the level. Starting with the irreducible
mod $\ell$ representation afforded by $F[\ell]$ of conductor either
$pq$ or $pqr$, and for any given $t \geq 4$,  we can find forms
giving the same residual representations in levels equal to the
product of $t$ primes by just raising the level $t-2$ ($t-3$,
respectively) times. For this, we have to take care not to lose
spurious primes as we add on new ones. The required analysis is
carried out in \cite{DT}.\\
 We conclude that Theorem \ref{squarefree} is true for any $t\geq 3$. Since
 the case of prime level was dealt with in the previous section, at
 this point only the case of $t=2$ remains unsolved (and a proof of
 Goldbach's conjecture would be enough to handle it).

\section{The case $t=2$ of Theorem \ref{squarefree} via a result of Ogg}
One way to treat Theorem \ref{squarefree} for the case where $N =
pq$ (without proving Goldbach's conjecture) is to appeal to the
results of Ogg in \cite{O}. If $p$ and $q$ are distinct primes, Ogg
finds a degree $0$ cuspidal divisor on $X_0(pq)$ whose image on
$J_0(pq)$ has order equal to the numerator of the fraction $(p -
1)(q + 1)/24$. Take $\ell > 3$. If $\ell$ divides $(p - 1)(q + 1)$,
using exactly the same arguments applied in (cf. \cite{Ma}) in the
case of prime level we find an eigenform $f$ at level $pq$ that is
Eisenstein mod $\ell$, therefore giving a reducible residual mod
$\ell$ representation. In particular, this means that the
coefficient $a_2$ satisfies again congruence (1) as in section 2,
for some prime $\lambda$ dividing $\ell$ in its field of
coefficients $K_f$, and we deduce as before that the degree of
$K_f$ is large (it grows with $\ell$).\\
We need to ensure that $f$ is genuinely a newform, i.e., that its
eigenvalues do not arise at level $p$ or at level $q$. We begin as
before by taking $\ell$ large. Then we find $q \equiv -1  \;
\mod{\ell}$ and pick $p$ to be a random prime that is not congruent
to $1$ mod $\ell$. Since the Eisenstein primes at prime level $N$
are divisors of $N -1$, we see that $\ell$ is not an Eisenstein
prime at either level $p$ or level $q$ while it is an Eisenstein
prime at level $pq$, thus the form $f$ must be a newform of level
$pq$. This completes the proof of the case $t=2$ of Theorem
\ref{squarefree}. Thus, putting together the results of the last
three sections, we conclude that Theorem \ref{squarefree} for any
positive value of $t$.


\section{The proof of Theorem \ref{almostsquarefree} }
To prove  Theorem \ref{almostsquarefree} we follow a strategy
similar to the one explained in section 3, except
 that now we will start from a diophantine equation such that the elliptic curve corresponding
  to any solution is a $Q$-curve defined over $\Q(i)$. For the Fermat-type equation $x^4+ y^2 = z^p$
  an attached $Q$-curve was proposed by Darmon in \cite{Da} and in the work of Ellenberg \cite{El} it
   was shown using the modularity of this curve that the diophantine equation does not have non-trivial
    solutions for large $p$. We will consider instead the diophantine problem:

$$ x^4 + y^2 =  5^{\ell} p $$
The result that we will prove in the next sections (see Theorem
\ref{main} in section 6, specialized to the case $c= 5^\ell$), which
is a generalization of the case $\ell = 0$ solved by
Friedlander-Iwaniec in \cite{fi2}, implies that for any prime
exponent $\ell$ there exist
infinitely many primes $p$ such that there are integer solutions $A,B$ to this equation. \\
Thus, if $\ell$ is a given prime and $A, B, p$ satisfy $$A^4 + B^2 = 5^\ell p $$
 with $p$ prime, we consider, as in the work of Darmon and Ellenberg, the elliptic curve $E$:
$$ y^2 = x^3 + 4 A x^2 + 2 (A^2 + i B) x $$
For simplicity, and since we have infinitely many primes $p$ satisfying the equation, we assume $p \neq \ell$.\\
The following properties of this curve are known (cf. \cite{El}):\\
$E$ is $2$-isogenous to its Galois conjugate, in particular it is a $Q$-curve. If we call $W$ its Weil restriction
 defined over $\Q$, it is a $\GL_2$-type abelian surface and thus it has a compatible family of $2$-dimensional
  Galois representations of $G_\Q$ attached. These representations have coefficients in $\Q(\sqrt{2})$.
   This abelian surface is semistable outside $2$ and it is modular. Computing the conductor of $W$ it
    follows that the modular form $f$ attached to $W$ has level $2^\alpha 5 p$, with $\alpha = 5$ or $8$.
     It has weight $2$ and trivial nebentypus. This newform has coefficients in $\Q(\sqrt{2})$ and it has
      an inner twist.
\\
We now consider for the prime $\ell$ we started with and $\lambda \mid \ell$ in $\Q(\sqrt{2})$ the
residual mod $\lambda$ representation $\bar{\rho}_{W,\lambda}$ attached to $W$. Assuming that $\ell > 13$
it is known that this representation is irreducible (cf. \cite{El}). Since the discriminant
of $E$ is $512 (A^2 + i B) 5^\ell p $ we can, as in \cite{El}, apply the Frey trick at
the semistable prime $5$ (observe that $5$ is unramified in $\Q(i) / \Q$): locally at $5$
 the valuation of the discriminant is divisible by $\ell$ (on the other hand, this does not
 happen locally at the prime $p$). Thus, we conclude that $\bar{\rho}_{W,\lambda}$ is unramified
 at $5$: more precisely it has conductor $2^\alpha p$ with $\alpha \in \{ 5 , 8\}$.
  If we apply lowering the level (cf. \cite{Ri}) we
 see that there is a newform $f'$ of level $2^\alpha p$, weight $2$ and trivial nebentypus such that some
 of its corresponding residual Galois representations in characteristic $\ell$ is $\bar{\rho}_{W,\lambda}$. \\
Now we can conclude as in section 3: If we look at the coefficient
$a_5$ of $f'$ since we know that there is a mod $\lambda$ congruence
with a newform $f$ corresponding to an abelian variety which is
semistable at $5$ (namely, the abelian surface $W$) then
 by the necessary condition for level raising we know that it must hold:
$$ a_ 5 \equiv \pm 6 \qquad \mod{\lambda} $$
From this congruence we see, as in sections 2 and 3, that the minimal field of definition of $a_5$, and a
 fortiori the field of coefficients of $f'$, has a degree that grows with $\ell$. This solves the case of
 almost-square-free level $2^\alpha p $ with $\alpha = 5$ or $8$, i.e., the case $t=1$, of Theorem \ref{almostsquarefree}.
  The  case of Theorem \ref{almostsquarefree} for levels of the
 form $2^\alpha p_1 \cdot ....\cdot p_t$, with $t > 1$ fixed and the $p_i$ odd, different primes, can be
 deduced from this by $t-1$ applications of raising the level as explained in section 3.

\section{Prime values of $(x^4 + y^2)/c$ }

We now introduce some notation which will be used from now on. For any given prime $p$, and any integer
$d$ we denote $v_p(d)$ the highest power of $p$ dividing $d$. Moreover,  we will write $d$ as
 $d=d_1d_2^2=d_1d_3^2d_4^4$, where $d_1,d_3$ are squarefree. We will consider
$\Lambda(r)$ the usual Von Mangoldt function, extended as zero over non integer numbers.
 Then, the main result of this section is the following theorem.
\begin{thm}\label{main} Let $B>0$ fixed. We have uniformly in $c\le (\log x)^B$
\begin{equation}
\sum\sum_{\kern-10pt(a^2+b^4)/c\le x}\Lambda((a^2+b^4)/c)=4\pi^{-1}\kappa G(c) (cx)^{3/4}+o\left(x^{3/4}\right)
\end{equation}
where $a,b$ run over positive integers, $G$ is a multiplicative function,
and
\begin{equation}\label{kappa}
\kappa=\int_0^1(1-t^4)^{1/2}dt=\Gamma(\frac{_{_1}}{^{^4}})^2/6\sqrt{2\pi}.
\end{equation}
\end{thm}
{\bf Remark:} It is important to emphasize here that the constant
$G(c)$, explicitly  described in Lemma \ref{multgeneral}, can take
the zero value, and it does  precisely in those $c$ which are non
representable as the sum of two squares, or such that $v_2(c)\equiv
3\tm 4$. For trivial reasons there is at most one prime in the
sequence $(a^2+b^4)/c$ in these cases, since none of the elements is
in fact coprime with $c$ if $c$ is non representable as the sum of
two squares, and if $v_2(c)\equiv 3\tm 4$, every integer of the
sequence has to be even. Hence, the proof that follows restricts to
those values of $c$ such that this constant is non zero since, in
any other case, the result is trivial.

\

Theorem \ref{main} is an easy generalization of Theorem 1.1 in  \cite{fi2}. However, some of
 the computations done in \cite{fi2} do not apply to this case in a straightforward manner and, hence,
   they must be done now with the required level of generality in the variable $c$.  In particular, the
   proof
of Theorem \ref{main} relies in the verification of the hypothesis needed to apply the Asymptotic Sieve
 due to Friedlander and Iwaniec in \cite{fi1}, but now for the sequence $a(c)_n=0$ for any $(n,c)>1$,
 and for $n$ coprime to $c$ given by
\begin{equation}\label{sequence}
\ac{c}_n=\sum\sum_{\kern-10pt(a^2+b^2)/c=n}\mathfrak Z(b),
\end{equation}
where $a,b$ are integers non necessarily positive, and $\mathfrak Z$
is the function with value $\mathfrak Z(m^2)=2$, for any integer
$m\ne 0$,  $\mathfrak Z(0)=1$, and $\mathfrak Z(b)=0$ in any other
case. From now on we will only consider integers $n$ coprime to $c$
and, then, we have  $\ac{c}_n= a^{old}_{cn}$ where $a^{old}_n$ is
the sequence related with Theorem 1 of \cite{fi2}. We now include
for reading convenience the hypotheses and main result of the
Asymptotic Sieve. The following, with the exception of
(\ref{divisor}), is basically a copy of Section 2 in \cite{fi2}. We
explain the difference between  (\ref{divisor}) and (2.8) in
\cite{fi2} at the end of this section.

\

Consider a sequence of real, nonnegative, numbers ${\cal A}=(a_n)_{n\ge 1}$, and $x$ a positive number.
We want to obtain an asymptotic formula for
$$
S(x)=\sum_{p\le x}a_p\log p,
$$
where the sum runs over prime numbers, in terms of $\A(x)=\sum_{n\le x} a_n$. We suppose
\begin{eqnarray}
\A(x)&\gg&\A(\sqrt{x})(\log x)^2\label{a:uno},\\
\A(x)&\gg&x^{1/3}\left(\sum_{n\le x}a_n^2\right)^{1/2}.\label{a:dos}
\end{eqnarray}
As usual in sieve theory, we will assume that for any integer $d>1$
$$
\A_d(x)=\sum_{{n\le x}\atop{d|n}}a_n=g(d)\A(x)+r_d(x),
$$
where $g$ is a multiplicative function, and $r_d(x)$ is regarded as an error term.
For the function $g$ we assume the following hypotheses
\begin{eqnarray}
&&0\le g(p^2)\le g(p)<1,\label{g:uno}\\
&&g(p)\ll p^{-1},\label{g:dos}\\
&&g(p^2)\ll p^{-2},\quad\text{ and }\label{g:tres}\\
&&\sum_{p\le y}g(p)=\log\log y+ e+O((\log y)^{-10}),\label{g:cuatro}
\end{eqnarray}
for every $y$ and some $e$ depending only on $g$. For the error term we will assume
\begin{equation}\label{super}
\sum_{d\le DL^2}^{\kern15pt 3}|r_d(t)|\le A(x)L^{-2},
\end{equation}
uniformly in $t\le x$, for some $D$ in the range
\begin{equation}\label{range:D}
 x^{2/3}<D<x.
\end{equation}
The superscript $3$ in (\ref{super}) restrict the sum to cube free moduli and $L=(\log x)^{2^{24}}$.
We also require
\begin{equation}\label{divisor}
 \A_d(x)\ll d^{-1}\tau(d)^8\A(x)\log x\text{ uniformly in }d\le x^{1/3},
\end{equation}
 and finally an estimate in bilinear forms like
\begin{equation}\label{bilinear}
 \sum_{m}|\sum_{{N<n\le 2N}\atop{{mn\le x}\atop{(n,m\Pi)=1}}}\beta(n)\mu(mn)a_{mn}|\le \A(x)L^{-4},
\end{equation}
where the coefficients are given by
\begin{equation}\label{coeff:bilinear}
\beta(n)=\beta(n,K)=\sum_{k|n,k\le K}\mu(k),
\end{equation}
for any $K$ in the range
\begin{equation}\label{range:K}
 1\le K\le xD^{-1},
\end{equation}
 $N$ verifing
\begin{equation}\label{range:N}
 \Delta^{-1}\sqrt{D}< N<\delta^{-1}\sqrt{x}
\end{equation}
for some $\Delta\ge\delta\ge 2$, and $\Pi$ is the product of all primes $p<P$ for some
$P$ which can be chosen conveniently in the range
\begin{equation}\label{range:P}
 2\le P\le\Delta^{1/2^{35}\log\log x}.
\end{equation}
In this conditions we have
\begin{prop}\label{prop1} Let $\cal A$ be a sequence verifying the above hypotheses. Then,
\begin{equation}
S(x)=H\A(x)\left\{1+O\left(\frac{\log \delta}{\log \Delta}\right)\right\}
\end{equation}
 where $H$ is the positive constant given by the convergent product
\begin{equation}\label{hache}
H=\prod_p(1-g(p))(1-\frac1p)^{-1},
\end{equation}
and the implied constant depends only on the function $g$.
\end{prop}
Normally $\delta$ is a large power of $\log x$ and $\Delta$ a small power of $x$.

{\bf Remark:} It is important to note that (\ref{divisor}) is not the original assumption (1.6)
in \cite{fi1}, but a slightly weaker. However, as the authors mention in that paper,
(1.6) is only required to reduce the hypotheses
\begin{equation}\label{bilinear:dos}
\sum_{m}\tau_5(m)|\sum_{{N<n\le 2N}\atop{{mn\le x}\atop{(n,m\Pi)=1}}}\beta(n)\mu(mn)a_{mn}|\le \A(x)(\log x)^{-3},
\end{equation}
 and
\begin{equation}\label{super:dos}
\sum_{d\le D}\mu^2(d)\tau_5(d)|r_d(t)|\le A(x)(\log x)^{-3},
\end{equation}
to (\ref{super}) and (\ref{bilinear}). We just have to follow the reasoning in Section 2, p. 1047 of
 \cite{fi1} to see that this reduction is also possible with our hypothesis (\ref{divisor}).
\section{Proof of Theorem \ref{main}}

To prove Theorem \ref{main} we will use Proposition \ref{prop1} for the sequence  given in (\ref{sequence}).
Hence, we have to check that the sequence verifies hypotheses (\ref{a:uno}) trhough (\ref{bilinear}).
Given an integer $d\ge 1$, we denote $\A_d(x;c)=\sum_{n\le x,n\equiv 0\tm d}\ac{c}_n$. The first thing
that  needs to be done is  to find a good aproximation of
$\A_d(x;c)$ in terms of a multiplicative function. Now,  $\A_d(x;c)=0$ for $(d,c)>1$ and for $(d,c)=1$ we have
\begin{equation}\label{moebius}
\A_d(x;c)=\sum_{k|c}\mu(k)\A^{old}_{ckd}(cx),
\end{equation}
and  we know by \cite{fi2} that
$$
\A^{old}_d(x)=g(d)\A^{old}(x)+r^{old}_d(x),
$$
where the functions $g$, $r^{old}$ satisfy conditions (\ref{a:uno}) trhough (\ref{bilinear}). Note
that for any integer $d$ the definition of $\A^{old}_{d}(x)$ is implicit in (\ref{moebius}) for $c=1$,
 and observe that $g$, $r^{old}$ are precisely the functions $g$, $r$ appearing in \cite{fi2}.  Hence,
  to approximate $\A_d(x;c)$ we  are tempted to   use the approximation of $\A^{old}(x)$ given in Lemma 3.4
  of \cite{fi2}. However, this lemma only works for cubefree integers $d$ which do not cover completely our case,
   since $c$ will be any number $c\le (\log x)^B$. Hence, our next objective is to generalize  Lemma 3.4 of \cite{fi2}
    to any integer $d$. As in \cite{fi2}, we start approximating $\A_d(x;c)$ by
$$
M_d(x;c)=\sum_{k|c}\mu(k)\frac1{ckd}\sum\sum_{\kern-10pt0<(a^2+b^2)\le cx}\mathfrak Z(b)\rho(b;ckd),
$$
for any $d$ coprime to $c$, where $\rho(b,d)$ denotes the number of solutions $\ap\tm d$ to the congruence
$\ap^2+b^2\equiv 0\tm{d}$,  and $M_d(x;c)=0$ otherwise. 
The following is a trivial consequence of Lemma 3.1 of \cite{fi2}.
\begin{lem}\label{error} Let $B>0$. For any $c\le (\log x)^B$ we have
$$
\sum_{d\le D}|\A_d(x;c)-M_d(x;c)|\ll D^{1/4}x^{9/16+\ep}
$$
for any $D\ge 1$ and $\ep>0$ and the implied constant depending only on $\ep$.
\end{lem}
Now, we need to find out the main term of $M_d(x;c)$, as we mentioned,  by generalizing Lemma 3.4 of \cite{fi2}.
\begin{lem}\label{multgeneral} Let $B>0$. We have uniformly for any $c\le (\log x)^B$
\begin{equation*}
M_d(x;c)=g_c(d)\left(4\kappa c^{3/4}G(c)\right) x^{3/4}+O\left(h(d)H(c)x^{1/2}\right),
\end{equation*}
where  $\kappa$ is given in (\ref{kappa}), $g_c(d)=0$ for any $(c,d)>1$, and  $g_c(d)=g(d)$
 otherwise where $g$ and  $h$ are the  multiplicative functions given by
\begin{equation}\label{defg}
g(d)=\frac1d\sum_{\nu_4|d_4}\nu_4^2\sum_{{\nu_3|d_3}\atop{(\nu_3,\frac{d_4}{\nu_4})=1}}\nu_3(\nu_3,d_1^*)\sum_{{\nu_1|\frac{d_1^*}{(\nu_3,d_1^*)}}\atop{(\nu_1,\frac{d_3d_4}{\nu_3\nu_4})=1}}
\rho\left(\delta(\nu_1,\nu_3,\nu_4)^2\right)\frac{\varphi\left(\delta(\nu_1,\nu_3,\nu_4)\right)}{d_1^*d_3d_4},
\end{equation}
and

\begin{equation}\label{defh}
\kern-2cm h(d)=\frac 1d\sum_{\nu_4|d_4}\nu_4^2\sum_{{\nu_3|d_3}\atop{(\nu_3,\frac{d_4}{\nu_4})=1}}\nu_3\sum_{{\nu_1|\frac{d_1^*}{(\nu_3,d_1^*)}}\atop{(\nu_1,\frac{d_3d_4}{\nu_3\nu_4})=1}}
\rho\left(\delta(\nu_1,\nu_3,\nu_4)^2\right)\tau(\delta(\nu_1,\nu_3,\nu_4))).
\end{equation}
Here $d_1^*=d_1/(d_1,2)$ and $\delta(\nu_1,\nu_3,\nu_4)=\frac{d_1^*d_3d_4}{(d_1^*,\nu_3)\nu_1\nu_3\nu_4}$. Finally
$G(c)=\sum_{k|c}\mu(k)g(ck)$ and  $H(c)=c^{1/2}\sum_{k|c}h(ck)$.  
\end{lem}
{\it Proof:}  We restrict only to integers $d$ coprime to $c$ since the result is trivial otherwise.
Given $d=d_1d_2^2=d_1d_3^2d_4^4$, with $d_1, d_3$ squarefree, and an integer $b$, let
us call $b_2=(b,d_2)$, $b_1=(b/b_2,d_1^*)$. Then, it is fairly straightforward to prove that
$$
\rho(b;d)=b_2\rho\left((d_1^*d_2/b_1b_2)^2\right)
$$
where $\rho$ is the multiplicative function given by
$$
\rho(p^\ap)=1+\chi_4(p),
$$
whit $\chi_4$ the character of conductor $4$, except $\rho(d)=0$ if $4|d$. Now, by definition, we have
\begin{equation}\label{summob}
M_d(x;c)=\sum_{k|c}\mu(k)M_{ckd}^{old}(cx)
\end{equation}
where
\begin{eqnarray*}
M_d^{old}(x)&=&\frac1{d}\sum\sum_{\kern-10pt0<(a^2+b^2)\le x}\mathfrak Z(b)\rho(b;d)=\frac2{d}\sum_{|r|\le x^{1/4}}\rho(r^2,d)\left\{ (x-r^4)^{1/2}+O(1)\right\}\\
&=&T_d+E_d,
\end{eqnarray*}
$T_d$  being the main term of $M_d^{old}(x)$, and
$$
E_d\ll\frac1{d}\sum_{|r|\le x^{1/4}}\rho(r^2,d)\ll x^{1/4}.
$$
Now suppose $(r,d_4)=\nu_4$. Then  $(r^2,d_3d_4^2)=\nu_4^2(r/\nu_4,d_3)$, and
 splitting the sum over the divisors of $d_4$ we get, after some calculations,
\begin{eqnarray*}
T_d(x)&=&\frac{2}{d}\sum_{\nu_4|d_4}\nu_4^4\sum_{{|r|\le z^{1/4}}\atop{(r,\frac{d_4}{\nu_4})=1}}(r,d_3)\rho\left(\left(\frac{d_1^*}{(r,d_1^*)}\frac{d_3}{(r,d_3)}\frac{d_4}{\nu_4}\right)^2\right)(z-r^4)^{1/2},
\end{eqnarray*}
where $z=z(\nu_4)=x/\nu_4^4$. Similarly, splitting the inner sum over the divisors of $d_3$ and  $d_1^*$,
it is easy to get
$$
T_d(x)=\frac{2}{d}\sum_{\nu_4|d_4}\nu_4^4\sum_{{\nu_3|d_3}\atop{(\nu_3,\frac{d_4}{\nu_4})=1}}
\nu_3^3\sum_{{\nu_1|\frac{d_1^*}{(\nu_3,d_1^*)}}\atop{(\nu_1,\frac{d_3d_4}{\nu_3\nu_4})=1}}
\nu_1^2
\rho\left(\delta(\nu_1,\nu_3,\nu_4)^2\right)
\sum_{{|r|\le u^{1/4}}\atop{{(r,\delta(\nu_1,\nu_3,\nu_4)=1}}}
(u-r^4)^{1/2},
$$
where $\delta(\nu_1,\nu_3,\nu_4)=\frac{d_1^*d_3d_4}{(\nu_3,d_1^*)\nu_1\nu_3\nu_4}$, and $u=u(\nu_1,\nu_3,\nu_4)=z/(\nu_3\nu_1)^4$. Estimating the inner sum, (also done in Lemma 3.4 of \cite{fi2}), we get
\begin{multline*}
T_d(x)=\frac{2}{d}\sum_{\nu_4|d_4}\nu_4^4\sum_{{\nu_3|d_3}\atop{(\nu_3,\frac{d_4}{\nu_4})=1}}\nu_3^3\sum_{{\nu_1|\frac{d_1^*}{(\nu_3,d_1^*)}}\atop{(\nu_1,\frac{d_3d_4}{\nu_3\nu_4})=1}}\nu_1^2
\rho\left(\delta(\nu_1,\nu_3,\nu_4)^2\right)\left\{\frac{\varphi\left(\delta(\nu_1,\nu_3,\nu_4)\right)}{\delta(\nu_1,\nu_3,\nu_4)}{2\kappa u^{3/4}}+\right.\\
\left. \frac{}{}O(\tau(\delta(\nu_1,\nu_3,\nu_4))u^{1/2})\right\}.
\end{multline*}
Substituting the value of $u$ in the previous formula  we get, for
any integer $d\ge 1$ not necessarily coprime with $c$,
$$
M_d^{old}(x)=g(d)4\kappa x^{3/4}+O(h(d)x^{1/2}),
$$
where $g$ and $h$ are the multiplicative functions given by (\ref{defg}) and (\ref{defh})
respectively. We just have to plug this into (\ref{summob}) to get the result with $g_c=g$
if $(d,c)=1$, and  $0$ otherwise. Note that
$$
G(c)=\sum_{k|c}\mu(k)g(ck)=\prod_{p|c}\left(g(p^{v_p(c)})-g(p^{v_p(c)+1})\right),
$$
defines a multiplicative function. Also,  observe that the
multiplicativity of both $g$ and $h$ is a direct consequence of the
definition in each case. It is straightforward to see that the value
at prime powers is given by $g(2^{4\ap+r})=1/2^{3\ap+r}$, and for
$p$ odd,
\begin{equation}\label{firstg}
g(p^{4\ap+r})=
\frac1{p^{4\ap+r}}(1+\chi_4(p))\left(1-\frac1{p}\right)\frac{p^{\ap}-1}{p-1}+\frac1{p^{3\ap+r}}g_r(p)
\end{equation}
where
\begin{equation}\label{secondg}
 g_{r}(p)=\begin{cases}1
&\text{ if } r=0\\
1+\chi_4(p)\left(1-\frac1p\right)&\text{ if } r=1\\
1+(1+\chi_4(p))\left(1-\frac1p\right)&\text{ if } r=2,3.\\
                \end{cases}
\end{equation}
For $h$ we will only need its value for cubefree integers which comes from the following
\begin{equation}\label{h:primes}
h(p)p=1+2\rho(p),\qquad h(p^2)p^2=p+2\rho(p),
\end{equation}
already gotten in Lemma 3.4 of \cite{fi2}.
%

\

We are now in position to verify hypotheses (\ref{a:uno}) trhough (\ref{bilinear}) for
the approximation
$$
 \A_d(x;c)=g_c(d)\A(x;c)+r_d(x).
$$
First of all we note that, by Lemmas \ref{error} and \ref{multgeneral},  we have
\begin{equation}\label{newa}
 \A_d(x;c)=g(d)4\kappa G(c) (cx)^{3/4}+O(d^{1/4}x^{9/16+\ep})+O(H(c)h(d)x^{1/2}),
\end{equation}
meanwhile, trivially, we have $a_c(n)\ll \tau(n)$. From here,  (\ref{a:uno}) and (\ref{a:dos}) follow
immediately by noting that $H(c)\ll c^\ep\ll \log x$, which is an easy consequence of
the definition of $h$. Also (\ref{g:uno}),(\ref{g:dos}),(\ref{g:tres}) and (\ref{g:cuatro}) are easy
 consequences of (\ref{firstg}),  (\ref{secondg}) and the Prime Number Theorem in the
arithmetic progression modulo $4$. Note  also that these conditions were already verified in
\cite{fi2} since, for cubefree integers, $g(d)$ is the same function as the one appearing in that
reference. In order to get (\ref{super}) we note that, by (\ref{newa}) used for any given $d$ and
 for $d=1$, we have
$$
r_d(x)=O(g(d)d^{1/4}x^{9/16+\ep})+O(H(c)h(d)x^{1/2}).
$$
Moreover
$$
\sum_{d\le x}^{\kern15pt 3}h(d)\ll \prod_{p\le x}(1+h(p))(1+h(p)^2)\ll (\log x)^6,
$$
and
$$
\sum_{d\le x}^{\kern15pt 3}g(d)\ll \prod_{p\le x}(1+g(p))(1+g(p)^2)\ll (\log x)^2,
$$
which gives
$$
\sum_{d\le D}^{\kern15pt 3}\left|r_d(t)\right|\ll D^{1/4}x^{9/16+\ep},
$$
and, in particular, implies (\ref{super}). For the remainder conditions, (\ref{divisor}) and
(\ref{bilinear}), we want to use the analogous results obtained in \cite{fi2} for the sequence
 $a^{old}_n$. It is then mandatory
to obtain the relation between the size of  $\A_d(x;c)$ and $\A^{old}_d(x)$. Now,
 (3.18) of \cite{fi2} is, with our notation
\begin{equation}\label{olda}
 \A^{old}(x)=4\kappa x^{3/4}+O(x^{1/2}).
\end{equation}
Hence, in view of (\ref{newa}), with $d=1$, to compare $\A^{old}(x)$ with $\A(x;c)$ we need
to control $G(c)$, and $G(c)^{-1}$ for any $c\le (\log x)^B$. This is the content of the next lemma.
\begin{lem}\label{Gc}  Let  $c$ be an integer. If $c$ is not representable as the sum of two squares,
 or $v_2(c)\equiv 3\tm 4$, then $G(c)=0$. Otherwise we have
$$
\frac1{c}\le G(c)\le \frac1{c^{3/4}}.
$$
\end{lem}
{\bf Remark:} The lemma is intended to show bounds which are enough for our purpose and by no
 means need to be optimal.

\

{\it Proof:} From  (\ref{firstg}) and  (\ref{secondg}) we see that
$$
G(p^{\ap})\ge \frac 1{p^\ap},
$$
which in particular implies the lower bound.
The upper bound is a direct consequence of (\ref{olda}) and (\ref{newa}) since, by definition,
 $\A(x;c)\le\A(x;1)=\A^{old}(x)$.

\

We now deal with (\ref{divisor}), for $(d,c)=1$. It is clear that in this case $\A_d(x;c)\le \A_{cd}^{old}(cx)$
 and it is trivial to get
$$
\A_d^{old}(x)\ll \frac1d\tau(d)\A^{old}(x),
$$
uniformly in $d\le x^{1/2-\ep}$. Hence, we have
$$
\A_d(x;c)\ll\frac1{cd}\tau(cd)\A^{old}(cx)\le\frac1{d}\tau(d)\tau(c)\A(x;c)\ll\frac1{d}\tau(d)\A(x;c)\log x,
$$
where we have used the lower bound in Lemma \ref{Gc}, together with (\ref{newa}) with $d=1$, (\ref{olda})
 and the bound $c\le(\log x)^B$. Hence, we are left with the bilinear condition
(\ref{bilinear}). We will  get this bound from the similar one
achieved in Proposition 4.1 of \cite{fi2}.

\


\begin{prop}\label{errorbilinear} Let $\eta>0$, $A>0$ and $B>0$. Then for any $c\le (\log x)^B$ we have
\begin{equation}\label{bilinear:concreta}
\sum_{m}|\sum_{{N<n\le 2N}\atop{{mn\le x}\atop{(n,m\Pi)=1}}}\mu(mn)\beta(n)\,\ac{c}_{mn}|\le \, \A(x;c)L^{4-A},
\end{equation}
for every $N$ with
\begin{equation}
 x^{1/4+\eta}<N<x^{1/2}(\log x)^{-U},
\end{equation}
where the coefficients $\beta(n)$ are given by (\ref{coeff:bilinear}) for any $1\le C\le N^{1-\eta}$,
$$
(\log\log x)^2\le \log P\le (\log x)(\log\log x)^{-2}.
$$
where $\Pi=\prod_{p<P}p$. Here $U$ and the implied constant in {\rm(\ref{bilinear:concreta})} need
to be taken sufficiently large in terms of $\eta$ and $A$.
\end{prop}

{\it Proof of Proposition \ref{errorbilinear}:}

\begin{eqnarray*}
&&\sum_{m}|\sum_{{N<n\le 2N}\atop{{mn\le x}\atop{(n,m\Pi)=1}}}\mu(mn)\beta(n)\,\ac{c}_{mn}|=
\sum_{(m,c)=1}|\sum_{{N<n\le 2N}\atop{{mn\le x}\atop{(n,m\Pi)=1}}}\mu(mn)\beta(n)\,a^{old}_{cmn}|\\
&\le&
\sum_{c|m}|\sum_{{N<n\le 2N}\atop{{mn\le cx}\atop{(n,m\Pi)=1}}}\mu(mn)\beta(n)\,a^{old}_{mn}|\le
\sum_{m}|\sum_{{N<n\le 2N}\atop{{mn\le cx}\atop{(n,m\Pi)=1}}}\mu(mn)\beta(n)\,a^{old}_{mn}|\\
&\le&
\A^{old}(x)L^{4-A},
\end{eqnarray*}
where we have used Proposition 4.1 of \cite{fi2}, and $c\le (\log x)^B$ together with the
lower bound on $P$. We just have to use (\ref{newa}) with $d=1$,  (\ref{olda}), and the lower bound
 in Lemma \ref{Gc} to get the result.

\

Theorem \ref{main} is now a direct consequence of Proposition \ref{prop1} and the upper bound
in Lemma \ref{Gc}. Recall that  $H$ is given by (\ref{hache}), which  in this case is
$$
H=\prod_p(1-\chi_4(p)p^{-1})=L(1,\chi_4)^{-1}=\frac 4\pi,
$$
and we are counting only positive integers $a,b$ in Theorem \ref{main}.

\noindent Luis V\'{i}ctor Dieulefait \newline Departament d'Algebra
i Geometria
\newline Universitat de Barcelona \newline
Gran Via de les Corts Catalanes 585 \newline 08007 Barcelona
\newline Spain
\newline email: ldieulefait@ub.edu
\\
\newline

\noindent Jorge Jim\'enez Urroz \newline Departament de Matem\`atica
Aplicada IV \newline Universitat Polit\`ecnica de Catalunya (UPC)
\newline Edifici C3-Campus Nord  \newline Jordi Girona, 1-3.
\newline E-08034 Barcelona \newline Spain \newline  email: jjimenez@ma4.upc.edu
\\
\newline

\noindent Kenneth Alan Ribet \newline Department of Mathematics 3840
\newline University of California
\newline Berkeley, CA 94720-3840 \newline USA \newline
email: ribet@math.berkeley.edu


\begin{thebibliography}{999}
\bibitem{Br} A. Brumer, The rank of $J_0(N)$.
Columbia University Number Theory Seminar (New York, 1992).
Ast\'erisque No. 228 (1995), 3, 41-68.

\bibitem{Ch} J. R. Chen,  On the representation of a larger even integer
as the sum of a prime and the product of at most two primes. Sci.
Sinica 16 (1973), 157-176.

\bibitem{Da} H. Darmon. Serre's conjectures. In V. Kumar Murty, editor, Seminar on Fermat's Last Theorem,
number 17 in CMS Conference Proceedings,  135-153, 1995.

\bibitem{DT} F. Diamond, R. Taylor,
Nonoptimal levels of mod $l$ modular representations. Invent. Math.
115 (1994), no. 3, 435-462.

\bibitem{El} J. S. Ellenberg,  Galois representations attached to $\mathbb Q$-curves and the
generalized Fermat equation $A^4+B^2=C^p$. Amer. J. Math. 126
(2004), no. 4, 763-787.

\bibitem{fi1} J.B. Friedlander, H. Iwaniec, Asymptotic sieve for primes, Ann. of Math. 148 (1998), 1041-1065.

\bibitem{fi2} J.B. Friedlander, H. Iwaniec, The polynomial $X^2 +Y^4$  captures its primes.
Ann. Math. 148 (1998), no. 3, 965-1040.

\bibitem{Ma} B. Mazur,  Modular curves and the Eisenstein ideal. Inst. Hautes \' Etudes Sci.
Publ. Math. No. 47 (1977), 33-186

\bibitem{O} A. P.  Ogg, Hyperelliptic modular curves. Bull. Soc. Math. France 102 (1974), 449-462.

\bibitem{Rilr} K. A. Ribet,  Raising the levels of modular representations. S\'{e}minaire de Th\'{e}orie des Nombres,
 Paris 1987-88, 259-271, Progr. Math., 81, Birkhauser Boston, Boston, MA, 1990.

\bibitem{Ri} K. A. Ribet,  On modular representations of ${\rm Gal}(\overline{\mathbb Q}/{\mathbb Q})$
 arising from modular forms. Invent. Math. 100 (1990), no. 2, 431-476.

 \bibitem{Sa} H. Saito, On a decomposition of spaces of cusp forms and trace formula of Hecke
operators. Nagoya Math. J. 80 (1980), 129-165.

\bibitem{Wi} A. Wiles, Modular elliptic curves and Fermat's last theorem.
Ann. of Math.  141 (1995), no. 3, 443-551.
\end{thebibliography}
\end{document}